\documentclass{amsart}
\usepackage[square]{natbib}
\newtheorem{thm}{Theorem}
\newtheorem{lem}[thm]{Lemma}

\theoremstyle{definition}

\newtheorem{ex}[thm]{Example}

\providecommand{\abs}[1]{\lvert#1\rvert}
\providecommand{\Abs}[1]{\Bigl\lvert#1\Bigr\rvert}

\begin{document}
%E PER PRIMA COSA SI METTONO TITOLO, AUTORI E VARIE DEL GENERE
\title[conditional 0-1 laws]{A conditional 0-1 law for the symmetric $\sigma$-field}
\author{Patrizia Berti}
\address{Patrizia Berti, Dipartimento di Matematica Pura ed Applicata ''G. Vitali'', Universita' di Modena e Reggio-Emilia, via Campi 213/B, 41100 Modena, Italy}
\email{berti.patrizia@unimore.it}

\author{Pietro Rigo}
\address{Pietro Rigo (corresponding author), Dipartimento di Economia Politica e Metodi Quantitativi, Universita' di Pavia, via S. Felice 5, 27100 Pavia, Italy}
\email{prigo@eco.unipv.it}

\begin{abstract} Let $(\Omega,\mathcal{B},P)$ be a
probability space, $\mathcal{A}\subset\mathcal{B}$ a
sub-$\sigma$-field, and $\mu$ a regular conditional distribution for
$P$ given $\mathcal{A}$. For various, classically interesting,
choices of $\mathcal{A}$ (including tail and symmetric) the
following 0-1 law is proved: There is a set $A_0\in\mathcal{A}$ such
that $P(A_0)=1$ and $\mu(\omega)(A)\in\{0,1\}$ for all
$A\in\mathcal{A}$ and $\omega\in A_0$. Provided $\mathcal{B}$ is
countably generated (and certain regular conditional distributions
exist), the result applies whatever $P$ is.
\end{abstract} \maketitle

\section{Introduction} \label{intro} Let $(\Omega,\mathcal{B},P)$ be a probability space,
$\mathcal{A}\subset\mathcal{B}$ a sub-$\sigma$-field, and
$\mathbb{P}$ the set of all probability measures on $\mathcal{B}$. A
{\em regular conditional distribution} ({\frenchspacing r.c.d.}) for
$P$ given $\mathcal{A}$ is a mapping
$\mu:\Omega\rightarrow\mathbb{P}$ such that $\mu(\cdot)(B)$ is a
version of $E(I_B\mid\mathcal{A})$ for all $B\in\mathcal{B}$.
Throughout, $P$ {\it is assumed to admit a {\frenchspacing r.c.d.}
given $\mathcal{A}$, denoted by $\mu$, and $\mathcal{B}$ is
countably generated} (that is, $\mathcal{B}$ is generated by one of
its countable subclasses).

\vspace{0.1cm}

We aim at showing that, for certain sub-$\sigma$-fields
$\mathcal{A}$ (including tail and symmetric), $\mu$ obeys the
following 0-1 law: There is a set $A_0\in\mathcal{A}$ with
$P(A_0)=1$ and
\begin{equation}\label{2law}
\mu(\omega)(A)\in\{0,1\}\quad\text{for all }A\in\mathcal{A}\text{
and }\omega\in A_0.
\end{equation}

\section{Motivations}
In the sequel, $A_0$ denotes a set of $\mathcal{A}$ satisfying
$P(A_0)=1$.

For both foundational and technical reasons, it would be desirable
that
\begin{equation}\label{uno}
\mu(\omega)(A)=I_A(\omega)\quad\text{for all } A\in\mathcal{A}\text{
and }\omega\in A_0
\end{equation}
for some $A_0$. Despite its heuristic content, however, condition
\eqref{uno} need not be true. In fact, by results of Blackwell and
Dubins (see \cite{BD} and references therein), condition \eqref{uno}
holds if and only if the trace $\sigma$-field $\mathcal{A}\cap
A_0=\{A\cap A_0:A\in\mathcal{A}\}$ is countably generated for some
$A_0$. Unless $\mathcal{A}$ is countably generated, thus,
\eqref{uno} does not hold for a number of probability measures $P$.

When \eqref{uno} fails, a natural question is whether some of its
consequences are still in force. The 0-1 law in \eqref{2law} is just
a (intriguing) consequence of condition \eqref{uno}.

To give \eqref{2law} some interpretation, let us fix
$\omega_0\in\Omega$. If $\mu(\omega_0)$ is 0-1 on $\mathcal{A}$,
then $\mu(\omega_0)(A\cap B)=\mu(\omega_0)(A)\mu(\omega_0)(B)$ for
$A\in\mathcal{A}$ and $B\in\mathcal{B}$. Conversely, the latter
relation (with $B=A$) yields $\mu(\omega_0)(A)=\mu(\omega_0)(A)^2$
for $A\in\mathcal{A}$, so that
\begin{equation*}
\mu(\omega_0)\text{ is 0-1 on
}\mathcal{A}\quad\Leftrightarrow\quad\mathcal{B}\text{ is
independent of }\mathcal{A}\text{ under }\mu(\omega_0).
\end{equation*}
Now, roughly speaking, the probability measure $\mu(\omega)$ should
embody the information conveyed by $\mathcal{A}$ for each $\omega$
in some set $A_0$. Accordingly, $\mathcal{B}$ should be independent
of $\mathcal{A}$, under $\mu(\omega)$, for all $\omega$ in such
$A_0$. This is precisely condition \eqref{2law}.

An equivalent (heuristic) argument is the following. If
$\mu(\omega_0)$ already includes the information in $\mathcal{A}$,
then $\mu$ should be a {\frenchspacing r.c.d.} for $\mu(\omega_0)$
given $\mathcal{A}$ as well. In fact, letting
$M=\{Q\in\mathbb{P}:\mu$ is a {\frenchspacing r.c.d.} for $Q$ given
$\mathcal{A}\}$, condition \eqref{2law} holds if and only if
$\mu(\omega)\in M$ for each $\omega$ in some $A_0$; see Theorem 12
of \cite{BR07}.

\begin{ex} \label{2exch} {\bf (Tail $\sigma$-field)} Suppose $\mathcal{A}=\bigcap_n\sigma(X_n,X_{n+1},\ldots)$ is the
tail $\sigma$-field of a sequence $(X_n)$ of real random variables
on $(\Omega,\mathcal{B},P)$. A probability measure $Q\in\mathbb{P}$
is 0-1 on $\mathcal{A}$ if and only if $\mathcal{B}$ is
asymptotically independent of $(X_n,X_{n+1},\ldots)$ under $Q$, in
the sense that
\begin{equation*}
\sup_{H\in\sigma(X_n,X_{n+1},\ldots)}\abs{Q(B\cap
H)-Q(B)Q(H)}\rightarrow 0\quad\text{for each }B\in\mathcal{B}.
\end{equation*}
Hence, condition \eqref{2law} becomes: $\mathcal{B}$ is
asymptotically independent of $(X_n,X_{n+1},\ldots)$, under
$\mu(\omega)$, for each $\omega$ in some $A_0$. This looks quite
reasonable (to us). Indeed, in \cite{BR07}, condition \eqref{2law}
is shown to be true (whatever $P$ is) if $\mathcal{A}$ is a tail
$\sigma$-field.
\end{ex}

A nice property of \eqref{2law} is that it is preserved under an
absolutely continuous change of probability measure.

\begin{ex}\label{abscont} {\bf (Absolute continuity)} Suppose
$Q\in\mathbb{P}$ satisfies $Q\ll P$ and $\mu(\omega)$ is 0-1 on
$\mathcal{A}$ for each $\omega$ in some $A_0$ (with
$A_0\in\mathcal{A}$ and $P(A_0)=1$). Let $f$ be a density of $Q$
with respect to $P$ and $A_1=\{\omega:0<\int
f(x)\mu(\omega)(dx)<\infty\}$. Then, $A_1\in\mathcal{A}$,
$Q(A_1)=1$, and
\begin{equation*}
\nu(\omega)(B)=\frac{\int_Bf(x)\mu(\omega)(dx)}{\int
f(x)\mu(\omega)(dx)},\quad B\in\mathcal{B},\omega\in A_1,
\end{equation*}
is a {\frenchspacing r.c.d.} for $Q$ given $\mathcal{A}$. If
$\omega\in A_0\cap A_1$, then $\nu(\omega)\ll\mu(\omega)$ and
$\mu(\omega)$ is 0-1 on $\mathcal{A}$. Thus,
$\nu(\omega)=\mu(\omega)$ on $\mathcal{A}$ for all $\omega\in
A_0\cap A_1$. This fact has two consequences. First, since
$Q(A_0\cap A_1)=1$, condition \eqref{2law} holds under $Q$ as well.
Second, if $P$ and $Q$ are equivalent (i.e. $Q\ll P$ and $P\ll Q$)
then $\nu=\mu$ on $\mathcal{A}$ {\frenchspacing a.s.}. This seems in
line with intuition.
\end{ex}

Since \eqref{2law} holds in various real situations, one could
suspect that it is always true, at least under mild conditions. This
is not so. Define in fact
\begin{equation*}
\mathcal{N}=\{B\in\mathcal{B}:P(B)=0\}.
\end{equation*}
Then, for \eqref{2law} to fail, it is enough that: (i)
$\mathcal{A}\supset\mathcal{N}$; (ii) $P\{x:\mu(x)=\mu(\omega)\}=0$
for each $\omega$ in some $A_0$; (iii) $P\{x:\mu(x)$ is not 0-1 on
$\mathcal{B}\}>0$; see Proposition 11 of \cite{BR07}. Conditions
(ii)-(iii) hold in most interesting problems. Thus, \eqref{2law}
typically fails whenever $\mathcal{A}\supset\mathcal{N}$. The next
two examples illustrate this fact.

\begin{ex}\label{frtyu} {\bf (A failure of condition \eqref{2law})}
Let $(\mathcal{F}_t:t\geq 0)$ be a filtration on
$(\Omega,\mathcal{B},P)$. As in stochastic calculus, suppose
$(\mathcal{F}_t)$ is right continuous and
$\mathcal{F}_0\supset\mathcal{N}$ (the so called "usual
conditions"). Suppose also that $\mathcal{B}$ is countably generated
and $X=\{X_t:t\geq 0\}$ is a real homogeneous Markov process,
relative to $(\mathcal{F}_t)$, with transition kernel
\begin{equation*}
K_t(a,H)=\text{Prob}\bigl(X_t\in H\mid X_0=a\bigr),\quad
t>0,a\in\mathbb{R},H\text{ a real Borel set}.
\end{equation*}
Letting $\mathcal{A}=\mathcal{F}_t$ for some $t>0$, one obtains
\begin{equation*}
\mu(\omega)\bigl(X_{2t}\in
\cdot\bigr)=K_t\bigl(X_t(\omega),\cdot\bigr)\quad\text{for each
}\omega\text{ in some set }A_0.
\end{equation*}
Thus, \eqref{2law} fails under various conditions on $K_t$. For
instance, \eqref{2law} fails whenever $K_t(a,\cdot)\neq
K_t(b,\cdot)$ for all $a\neq b$ and $K_t(a,\{b\})=0$ for all $a,b$.
In fact, $\mu(\omega)$ is not 0-1 on $\mathcal{B}$ for each
$\omega\in A_0$. Moreover, $P(X_t=b)=\int K_t(X_0,\{b\})dP=0$ for
all $b$. Hence,
\begin{equation*}
P\{x:\mu(x)=\mu(\omega)\}\leq
P\bigl(X_t=X_t(\omega)\bigr)=0\quad\text{for all }\omega\in A_0.
\end{equation*}
Therefore, conditions (i)-(ii)-(iii) hold.
\end{ex}

Incidentally, Example \ref{frtyu} also suggests the following remark
(unrelated to condition \eqref{2law}). According to a usual naive
interpretation, $\mathcal{F}_t$ describes the information at time
$t$, in the sense that each event $A\in \mathcal{F}_t$ is known to
be true or false at time $t$. This interpretation does not make
sense in Example \ref{frtyu} as far as $\{X=x\}\in\mathcal{B}$ for
each possible path $x$ of the process $X$. In fact, $P(X=x)\leq
P(X_t=x(t))=0$ for all $x$, so that $\{X=x\}\in\mathcal{F}_0$ for
every path $x$. Under such interpretation, thus, the X-path would be
already known at time $t=0$.

\begin{ex}\label{agmn} {\bf (One more failure of condition \eqref{2law}; see \cite{BR07})} Let $\Omega=\mathbb{R}^2$, $\mathcal{B}$ the Borel $\sigma$-field,
and $P=Q\times Q$ where $Q$ is the $N(0,1)$ law on the real Borel
sets. Define
$\mathcal{A}=\sigma\bigr(\mathcal{G}\cup\mathcal{N}\bigl)$ where
$\mathcal{G}$ is the $\sigma$-field on $\Omega$ generated by
$(x,y)\mapsto x$. A {\frenchspacing r.c.d.} for $P$ given
$\mathcal{G}$ is $\mu((x,y))=\delta_x\times Q$. Since
$\mathcal{A}=\sigma\bigr(\mathcal{G}\cup\mathcal{N}\bigl)$, $\mu$ is
also a {\frenchspacing r.c.d.} for $P$ given $\mathcal{A}$.
Moreover, for all $(x,y)$, one has $\{x\}\times
[0,\infty)\in\mathcal{A}$ and
\begin{equation*}
\mu((x,y))\Bigl(\{x\}\times [0,\infty)\Bigr)=\frac{1}{2}.
\end{equation*}
\end{ex}

Though implicit in ideas of Dynkin \cite{D} and Diaconis and
Freedman \cite{DF81}, condition \eqref{2law} has been almost
neglected so far. Possible related references are \cite{BR1999},
\cite{BR07}, \cite{BD} and \cite{SSK}, but only \cite{BR07} is
explicitly devoted to \eqref{2law}.

This note carries on the investigation started in \cite{BR07}. It is
proved that condition \eqref{2law} holds (whatever $P$ is) for
certain sub-$\sigma$-fields $\mathcal{A}$, including the symmetric
one.

\section{Results} Let $F$ be a class of measurable
functions $f:\Omega\rightarrow\Omega$, where measurability means
$f^{-1}\bigl(\mathcal{B}\bigr)\subset\mathcal{B}$. In case $F$ is a
group under composition, with the identity map on $\Omega$ as
group-identity, we briefly say that $F$ {\it is a group}. Whether or
not $F$ is a group, the $F$-{\it invariant} $\sigma$-field is
\begin{equation*}
\mathcal{A}_F=\{B\in\mathcal{B}:f^{-1}B=B\text{ for all }f\in F\}
\end{equation*}
and a probability measure $Q$ on $\mathcal{B}$ is $F$-invariant if
$Q\circ f^{-1}=Q$ for all $f\in F$. Let $\mathbb{P}_F$ denote the
set of $F$-invariant probability measures.

One more definition is to be recalled. Let
$\mathbb{Q}\subset\mathbb{P}$ be a collection of probability
measures and $\mathcal{G}\subset\mathcal{B}$ a sub-$\sigma$-field.
Then, $\mathcal{G}$ is {\em sufficient for} $\mathbb{Q}$ in case,
for each $B\in\mathcal{B}$, there is a $\mathcal{G}$-measurable
function $h:\Omega\rightarrow\mathbb{R}$ which is a version of
$E_Q\bigl(I_B\mid\mathcal{G}\bigr)$ for all $Q\in\mathbb{Q}$. When
$\mathbb{Q}=\mathbb{P}_F$, sufficiency is a key ingredient for
integral representation of invariant measures; see \cite{DF81},
\cite{D}, \cite{F} and \cite {M}. Conditions under which
$\mathcal{A}_F$ is sufficient for $\mathbb{P}_F$ are given in
Theorem 3 of \cite{F}. In particular, $\mathcal{A}_F$ is sufficient
for $\mathbb{P}_F$ if $F$ is a countable group or if $F$ includes
only one function.

Arguing as Maitra in \cite{M}, we now prove that condition
\eqref{2law} holds whenever $P\in\mathbb{P}_F$, $F$ is countable and
$\mathcal{A}$ sufficient for $\mathbb{P}_F$. Basing on this fact we
subsequently show that, if $F$ is a finite group and
$\mathcal{A}=\mathcal{A}_F$, then \eqref{2law} holds whatever $P$
is.
\begin{lem}\label{maitra}
Suppose $\mathcal{B}$ is countably generated, $F$ is countable, and
$P\ll P_0$ for some $P_0\in\mathbb{P}_F$ which admits a
{\frenchspacing r.c.d.} given $\mathcal{A}$. Then, condition
\eqref{2law} holds provided $\mathcal{A}$ is sufficient for
$\mathbb{P}_F$. In particular, when $\mathcal{A}=\mathcal{A}_F$,
condition \eqref{2law} holds if $F$ is a group or if $F$ includes
only one function.
\end{lem}
\begin{proof}
By Example \ref{abscont}, it is enough to prove that $P_0$ meets
\eqref{2law}. Thus, it can be assumed $P\in\mathbb{P}_F$. Let
$\mathcal{M}=\{B\in\mathcal{B}:Q(B)=0$ for all $Q\in\mathbb{P}_F\}$.
Since $\mathcal{A}$ is sufficient for $\mathbb{P}_F$ and
$\mathcal{B}$ is countably generated, by Theorem 1 of \cite{B},
there is a countably generated $\sigma$-field $\mathcal{D}$ such
that
$\mathcal{D}\subset\mathcal{A}\subset\sigma\bigl(\mathcal{D}\cup\mathcal{M})$.
Since $\mathcal{D}$ is countably generated and
$\mathcal{D}\subset\mathcal{A}$, there is $A_1\in\mathcal{A}$ such
that $P(A_1)=1$ and $\mu(\omega)(D)=I_D(\omega)$ for all $D
\in\mathcal{D}$ and $\omega\in A_1$. Since $P\in\mathbb{P}_F$, given
$B\in\mathcal{B}$ and $f\in F$, one obtains
$\mu(\omega)(f^{-1}B)=\mu(\omega)(B)$ for almost all $\omega$. Since
$F$ is countable and $\mathcal{B}$ countably generated, it follows
that $\mu(\omega)\in\mathbb{P}_F$ for each $\omega$ in some set
$A_2\in\mathcal{A}$ with $P(A_2)=1$. Fix $\omega\in A_1\cap A_2$.
Then, $\mu(\omega)$ is 0-1 on $\mathcal{D}\cup\mathcal{M}$, so that
$\mu(\omega)$ is 0-1 on $\sigma\bigl(\mathcal{D}\cup\mathcal{M})$ as
well. Since
$\mathcal{A}\subset\sigma\bigl(\mathcal{D}\cup\mathcal{M})$, for
getting condition \eqref{2law} it suffices to let $A_0=A_1\cap A_2$.
\end{proof}

\begin{thm}\label{group} If $\mathcal{B}$ is countably
generated, $F$ is a finite group and $\mathcal{A}=\mathcal{A}_F$,
then condition \eqref{2law} holds.
\end{thm}
\begin{proof} Define
\begin{equation*}
Q=\frac{\sum_{f\in F}P\circ
f^{-1}}{card(F)},\quad\nu(\omega)=\frac{\sum_{f\in
F}\mu(\omega)\circ f^{-1}}{card(F)}\quad\text{for all
}\omega\in\Omega,
\end{equation*}
and note that $Q=P$ and $\nu(\omega)=\mu(\omega)$ on
$\mathcal{A}=\mathcal{A}_F$. Further, $\nu(\omega)\in\mathbb{P}$ for
all $\omega\in\Omega$, $\omega\mapsto\nu(\omega)(B)$ is
$\mathcal{A}$-measurable for all $B\in\mathcal{B}$, and for each
$A\in\mathcal{A}$ and $B\in\mathcal{B}$ one obtains:
\begin{gather*}
\int_A\nu(\omega)(B)Q(d\omega)=\int_A\nu(\omega)(B)P(d\omega)\\=\frac{1}{card(F)}\sum_{f\in
F}\int_A\mu(\omega)(f^{-1}B)P(d\omega)=\frac{1}{card(F)}\sum_{f\in
F}P\bigl(A\cap f^{-1}B\bigr)\\=\frac{1}{card(F)}\sum_{f\in
F}P\bigl(f^{-1}A\cap f^{-1}B\bigr)=Q(A\cap B).
\end{gather*}
Hence, $\nu$ is a {\frenchspacing r.c.d.} for $Q$ given
$\mathcal{A}$. Since $F$ is a finite group, then $Q\in\mathbb{P}_F$
and $\mathcal{A}=\mathcal{A}_F$ is sufficient for $\mathbb{P}_F$.
Accordingly, by applying Lemma \ref{maitra} to $Q$ and $\nu$, there
is $A_0\in\mathcal{A}$ such that $Q(A_0)=1$ and $\nu(\omega)$ is 0-1
on $\mathcal{A}$ for each $\omega\in A_0$. Since $Q=P$ and $\nu=\mu$
on $\mathcal{A}$, this concludes the proof.
\end{proof}

Among other things, given a single measurable function
$f:\Omega\rightarrow\Omega$, Lemma \ref{maitra} and Theorem
\ref{group} apply to
$\mathcal{A}=\mathcal{A}_{\{f\}}=\{B\in\mathcal{B}:f^{-1}B=B\}$.
Precisely, Lemma \ref{maitra} grants condition \eqref{2law} in case
$P\ll P_0$ for some $f$-invariant $P_0$ (admitting a {\frenchspacing
r.c.d.} given $\mathcal{A}$). By Theorem \ref{group}, instead,
condition \eqref{2law} holds whatever $P$ is in case $f$ is
bijective with $f=f^{-1}$.

Towards our main example, concerning the symmetric $\sigma$-field
(cf. Example \ref{gxmmj}), we mention one more consequence of
Theorem \ref{group}.

\begin{ex}\label{xsetrr} {\bf (Permutations of order $n$)} Fix a measurable space $(\mathcal{X},\mathcal{U})$,
with $\mathcal{U}$ countably generated, and define
$(\Omega,\mathcal{B})=(\mathcal{X}^\infty,\mathcal{U}^\infty)$.
Denote points of $\Omega=\mathcal{X}^\infty$ by
$\omega=(\omega_1,\omega_2,\ldots)$. A {\it permutation of order}
$n$ is a map $f:\Omega\rightarrow\Omega$ of the form
\begin{equation*}
f(\omega)=(\omega_{\pi_1},\ldots,\omega_{\pi_n},\omega_{n+1},\ldots),\quad\omega\in\Omega,
\end{equation*}
for some permutation $(\pi_1,\ldots,\pi_n)$ of $(1,\ldots,n)$. The
set $F_n$ of permutations of order $n$ is a group with $n!$
elements, and $\mathcal{A}_{F_n}$ includes those $B\in\mathcal{B}$
invariant under permutations of the first $n$ coordinates. By
Theorem \ref{group}, every {\frenchspacing r.c.d.} $\mu$ (for some
law $P\in\mathbb{P}$) given $\mathcal{A}=\mathcal{A}_{F_n}$ meets
condition \eqref{2law}.
\end{ex}

We now turn to our main result. Let
$\mathcal{A}_n\subset\mathcal{B}$ be a sub-$\sigma$-field,
$n=1,2,\ldots$, and
\begin{equation*}\mathcal{A}_*=\sigma\bigl(\bigcup_{n\geq
1}\bigcap_{j\geq
n}\mathcal{A}_j\bigr),\quad\mathcal{A}^*=\sigma\bigl(\bigcap_{n\geq
1}\bigcup_{j\geq n}\mathcal{A}_j\bigr).
\end{equation*}
It seems reasonable that condition \eqref{2law} holds provided it
holds for every $\mathcal{A}_n$ and
$\mathcal{A}_n\rightarrow\mathcal{A}$ in some sense. In fact, this
is true if $\mathcal{A}\subset\mathcal{A}_*$ and
$\mathcal{A}_n\rightarrow\mathcal{A}$ is meant as
\begin{equation}\label{albereta}
E\bigl(I_B\mid\mathcal{A}_n\bigr)\overset{P}\rightarrow
E\bigl(I_B\mid\mathcal{A}\bigr)\quad\text{for each }B\in\mathcal{B}.
\end{equation}
Furthermore, condition \eqref{2law} holds for
$\mathcal{A}\subset\mathcal{A}^*$ (and not only for
$\mathcal{A}\subset\mathcal{A}_*$) if \eqref{albereta} is
strengthened into
\begin{equation}\label{reds}
E\bigl(I_B\mid\mathcal{A}_n\bigr)\overset{a.s.}\rightarrow
E\bigl(I_B\mid\mathcal{A}\bigr)\quad\text{for each
}B\in\mathcal{B}.\tag{3*}
\end{equation}
Note that, by the martingale convergence theorem, if $\mathcal{A}_n$
is a monotonic sequence then \eqref{reds} holds with
$\mathcal{A}=\mathcal{A}_*=\mathcal{A}^*$.

\begin{thm}\label{mwuv} Suppose $\mathcal{B}$ is countably generated
and, for each $n\geq 1$:
\begin{gather}\label{albergo}
\text{There are a {\frenchspacing r.c.d.} }\nu_n\text{ for }P\text{
given }\mathcal{A}_n\text{ and a set
}C_n\in\mathcal{A}_n\\\text{such that }P(C_n)=1\text{ and
}\nu_n(\omega)\text{ is 0-1 on }\mathcal{A}_n\text{ for all
}\omega\in C_n. \notag
\end{gather}
If \eqref{albereta} holds and $\mathcal{A}\subset\mathcal{A}_*$, or
if \eqref{reds} holds and $\mathcal{A}\subset\mathcal{A}^*$, then
\begin{equation*}
\mu(\omega)\text{ is 0-1 on }\mathcal{A}\text{ for each }\omega\in
A_0,\text{ where }A_0\in\mathcal{A}\text{ and }P(A_0)=1
\end{equation*}
(that is, condition \eqref{2law} holds). In particular, condition
\eqref{2law} holds whenever $\mathcal{A}_n$ is a monotonic sequence
and $\mathcal{A}=\mathcal{A}_*=\mathcal{A}^*$.
\end{thm}
\begin{proof} Suppose $\mathcal{A}\subset\mathcal{A}_*$ and \eqref{albereta} holds. Define
\begin{equation*}V_n^B(\omega)=\sup_{H\in\mathcal{A}_n}\Abs{\mu(\omega)(B\cap
H)-\mu(\omega)(B)\mu(\omega)(H)},\quad n\geq
1,B\in\mathcal{B},\omega\in\Omega,
\end{equation*}
and let $\mathcal{B}_0$ be a countable field such that
$\mathcal{B}=\sigma(\mathcal{B}_0)$. It is enough proving that:
\begin{gather}\label{citrew}
\text{There are a subsequence }(n_j)\text{ and a set
}A_0\in\mathcal{A}\text{ such that }\\P(A_0)=1\text{ and
}\lim_jV_{n_j}^B(\omega)=0\text{ for all }\omega\in A_0\text{ and
}B\in\mathcal{B}_0\notag.
\end{gather}

Suppose in fact \eqref{citrew} holds. Fix $\omega\in A_0$,
$B\in\mathcal{B}$ and $\epsilon>0$. Since $\mathcal{B}_0$ is a field
which generates $\mathcal{B}$, there is $B_0\in\mathcal{B}_0$ such
that $\mu(\omega)(B\Delta B_0)<\epsilon$. Hence,
\begin{gather*}
V_n^B(\omega)\leq\sup_{H\in\mathcal{A}_n}\Abs{\mu(\omega)(B\cap
H)-\mu(\omega)(B_0\cap
H)}+\\+\sup_{H\in\mathcal{A}_n}\Abs{\mu(\omega)(B_0\cap
H)-\mu(\omega)(B_0)\mu(\omega)(H)}+\sup_{H\in\mathcal{A}_n}\Abs{\mu(\omega)(B_0)\mu(\omega)(H)
-\mu(\omega)(B)\mu(\omega)(H)}\\\leq
V_n^{B_0}(\omega)+2\mu(\omega)(B\Delta
B_0)<V_n^{B_0}(\omega)+2\epsilon\quad\text{for all }n.
\end{gather*}
Since $\omega\in A_0$ and $B_0\in\mathcal{B}_0$, condition
\eqref{citrew} yields
\begin{equation*} \limsup_j V_{n_j}^B(\omega)\leq
2\epsilon+\limsup_j V_{n_j}^{B_0}(\omega)=2\epsilon.
\end{equation*}
Thus, $\lim_jV_{n_j}^B(\omega)=0$ for all $B\in\mathcal{B}$ and
$\omega\in A_0$. Denote $\mathcal{C}=\bigcup_{n\geq 1}\bigcap_{j\geq
n}\mathcal{A}_j$ and fix $\omega\in A_0$ and $A\in\mathcal{C}$. If
$A\in\mathcal{A}_k$ for some $k$, then
\begin{equation*}\abs{\mu(\omega)(A)-\mu(\omega)(A)^2}\leq\sup_{H\in\mathcal{A}_k}\Abs{\mu(\omega)(A\cap
H)-\mu(\omega)(A)\mu(\omega)(H)}=V_k^A(\omega).
\end{equation*}
Since $A\in\mathcal{C}$, there is $n$ such that $A\in\mathcal{A}_k$
for each $k\geq n$, so that
\begin{equation*}\abs{\mu(\omega)(A)-\mu(\omega)(A)^2}\leq\lim_jV_{n_j}^A(\omega)=0.
\end{equation*}
Therefore, $\mu(\omega)$ is 0-1 on $\mathcal{C}$, which implies that
$\mu(\omega)$ is 0-1 on $\sigma(\mathcal{C})=\mathcal{A}_*$. Since
$\mathcal{A}\subset\mathcal{A}_*$, condition \eqref{2law} holds.

It remains to prove condition \eqref{citrew}. The proof is split
into three steps.

\vspace{0.15cm}

{\bf (i)} Fix $n$ and take $\nu_n$ and $C_n$ as in condition
\eqref{albergo}. Since $C_n\in\mathcal{A}_n$ and $P(C_n)=1$, up to
modifying $\nu_n$ on $C_n^c$, it can be assumed that $\nu_n(\omega)$
is 0-1 on $\mathcal{A}_n$ for all $\omega\in\Omega$. We now prove
that, for each $\omega$ in some set $M_n\in\mathcal{A}$ with
$P(M_n)=1$, one has
\begin{equation*}
\mu(\omega)(B\cap H)=\int_{\{\nu_n(H)=1\}}
\nu_n(x)(B)\mu(\omega)(dx)\quad\text{for all
}H\in\mathcal{A}_n\text{ and }B\in\mathcal{B}
\end{equation*}
where $\{\nu_n(H)=1\}$ denotes the set $\{x:\nu_n(x)(H)=1\}$. Define
\begin{equation*}
\mu_n(\omega)(B)=\int\nu_n(x)(B)\mu(\omega)(dx),\quad\omega\in\Omega,
B\in\mathcal{B}.
\end{equation*}
Since $\mathcal{B}$ is countably generated and $\mu_n$ is a
{\frenchspacing r.c.d.} for $P$ given $\mathcal{A}$, there is
$M_n\in\mathcal{A}$ such that $P(M_n)=1$ and
$\mu_n(\omega)=\mu(\omega)$ for all $\omega\in M_n$. Let
$H\in\mathcal{A}_n$, $B\in\mathcal{B}$ and $\omega\in M_n$. Since
$\nu_n(\cdot)(H)\in\{0,1\}$, then $\nu_n(x)(B\cap
H)=\nu_n(x)(B)I_{\{\nu_n(H)=1\}}(x)$ for all $x\in\Omega$. Thus,
\begin{gather*}
\mu(\omega)(B\cap H)=\mu_n(\omega)(B\cap H)=\int\nu_n(x)(B\cap
H)\mu(\omega)(dx)\\=\int_{\{\nu_n(H)=1\}}\nu_n(x)(B)\mu(\omega)(dx).
\end{gather*}

\vspace{0.15cm}

{\bf (ii)} We next prove that, for each $\omega\in M_n\cap T$, where
$T\in\mathcal{A}$ and $P(T)=1$, one also has
\begin{equation*}
\mu(\omega)(B)\mu(\omega)(H)=\int_{\{\nu_n(H)=1\}}
\mu(x)(B)\mu(\omega)(dx)\quad\text{for all }H\in\mathcal{A}_n\text{
and }B\in\mathcal{B}.
\end{equation*}
Let $\sigma(\mu)$ be the $\sigma$-field generated by $\mu(\cdot)(B)$
for all $B\in\mathcal{B}$. Then, $\sigma(\mu)\subset\mathcal{A}$ and
$\sigma(\mu)$ is countably generated since $\mathcal{B}$ is
countably generated. Hence, there is $T\in\mathcal{A}$ with $P(T)=1$
and $\mu(\omega)(D)=I_D(\omega)$ for all $D\in\sigma(\mu)$ and
$\omega\in T$. Given $C\in\mathcal{B}$ and a bounded
$\sigma(\mu)$-measurable function $h:\Omega\rightarrow\mathbb{R}$,
it follows that
\begin{equation*}
\int_C
h(x)\mu(\omega)(dx)=h(\omega)\mu(\omega)(C)\quad\text{whenever
}\omega\in T.
\end{equation*}
Fix $\omega\in M_n\cap T$, $H\in\mathcal{A}_n$ and
$B\in\mathcal{B}$. Letting $h(x)=\mu(x)(B)$ and $C=\{\nu_n(H)=1\}$,
one obtains
\begin{gather*}
\mu(\omega)(B)\mu(\omega)(H)=\mu(\omega)(B)\mu(\omega)\bigl(\nu_n(H)=1\bigr)\quad\text{since
}\omega\in M_n\\=\int_{\{\nu_n(H)=1\}}
\mu(x)(B)\mu(\omega)(dx)\quad\text{since }\omega\in T.
\end{gather*}

\vspace{0.15cm}

{\bf (iii)} Since $\mathcal{B}_0$ is countable, by \eqref{albereta}
and a diagonalization argument, there is a subsequence $(n_j)$ such
that
\begin{equation*}
E\bigl(I_B\mid\mathcal{A}_{n_j}\bigr)\overset{a.s.}\rightarrow
E\bigl(I_B\mid\mathcal{A}\bigr),\text{ as }j\rightarrow\infty,\text{
for all }B\in\mathcal{B}_0.
\end{equation*} Define $Z_n^B(\cdot)=\nu_n(\cdot)(B)-\mu(\cdot)(B)$ .
If $B\in\mathcal{B}_0$, then $Z_{n_j}^B\overset{a.s.}\rightarrow 0$
as $j\rightarrow\infty$, and $\abs{Z_{n_j}^B}\leq 1$ for all $j$.
Hence,
\begin{equation*}
\int\abs{Z_{n_j}^B(x)}\mu(\cdot)(dx)=E\bigl(\abs{Z_{n_j}^B}\mid\mathcal{A}\bigr)\overset{a.s.}\rightarrow
0.
\end{equation*}
Define further
\begin{equation*}
S=\{\omega:\lim_j\int\abs{Z_{n_j}^B(x)}\mu(\omega)(dx)=0\text{ for
each }B\in\mathcal{B}_0\}\quad\text{and}\quad A_0=\bigcap_n (M_n\cap
S\cap T).
\end{equation*}
Then, $A_0\in\mathcal{A}$ and $P(A_0)=1$. Given $B\in\mathcal{B}_0$
and $\omega\in A_0$, points (i)-(ii) yield:
\begin{gather*}
V_{n_j}^B(\omega)=\sup_{H\in\mathcal{A}_{n_j}}\Abs{\mu(\omega)(B\cap
H)-\mu(\omega)(B)\mu(\omega)(H)}\\=\sup_{H\in\mathcal{A}_{n_j}}\Abs{\int_{\{\nu_{n_j}(H)=1\}}
\nu_{n_j}(x)(B)\mu(\omega)(dx)-\int_{\{\nu_{n_j}(H)=1\}}
\mu(x)(B)\mu(\omega)(dx)}\\=\sup_{H\in\mathcal{A}_{n_j}}\Abs{\int_{\{\nu_{n_j}(H)=1\}}
Z_{n_j}^B(x)\mu(\omega)(dx)}\\\leq\int
\abs{Z_{n_j}^B(x)}\mu(\omega)(dx)\rightarrow 0\quad\text{as
}j\rightarrow\infty.
\end{gather*}
Thus \eqref{citrew} holds, and this concludes the proof in case
$\mathcal{A}\subset\mathcal{A}_*$ and \eqref{albereta} holds.

Finally, suppose $\mathcal{A}\subset\mathcal{A}^*$ and \eqref{reds}
holds. By using \eqref{reds} instead of \eqref{albereta}, in point
(iii) there is no need of taking a subsequence $(n_j)$, and one
obtains $\lim_nV_n^B(\omega)=0$ for all $B\in\mathcal{B}_0$ and
$\omega$ in a set $A_0\in\mathcal{A}$ with $P(A_0)=1$. Arguing as at
the beginning of this proof, this in turn implies
$\lim_nV_n^B(\omega)=0$ for all $B\in\mathcal{B}$ and $\omega\in
A_0$. Denote $\mathcal{L}=\bigcap_{n\geq 1}\bigcup_{j\geq
n}\mathcal{A}_j$ and fix $\omega\in A_0$ and $A\in\mathcal{L}$.
Since $A\in\mathcal{L}$, there is a subsequence $(m_j)$ (possibly
depending on $A$) such that $A\in\mathcal{A}_{m_j}$ for all $j$.
Hence,
\begin{equation*}
\abs{\mu(\omega)(A)-\mu(\omega)(A)^2}\leq\sup_{H\in\mathcal{A}_{m_j}}\Abs{\mu(\omega)(A\cap
H)-\mu(\omega)(A)\mu(\omega)(H)}=V_{m_j}^A(\omega)\rightarrow 0.
\end{equation*}
Therefore, $\mu(\omega)$ is 0-1 on $\mathcal{L}$, which implies that
$\mu(\omega)$ is 0-1 on $\sigma(\mathcal{L})=\mathcal{A}^*$. Since
$\mathcal{A}\subset\mathcal{A}^*$, this concludes the proof.
\end{proof}

\vspace{0.3cm}

For Theorem \ref{mwuv} to apply, it is useful to have some condition
implying existence of {\frenchspacing r.c.d.'s} whatever the
sub-$\sigma$-field is. One such condition is: $P$ {\em admits a
{\frenchspacing r.c.d.} given $\mathcal{G}$, for any
sub-$\sigma$-field $\mathcal{G}\subset\mathcal{B}$, provided $P$ is
perfect and $\mathcal{B}$ countably generated}; see \cite{J}. We
recall that $P$ is {\em perfect} in case each
$\mathcal{B}$-measurable function $f:\Omega\rightarrow\mathbb{R}$
meets $P(f\in I)=1$ for some real Borel set $I\subset f(\Omega)$.
For $P$ to be perfect, it is enough that $\Omega$ is a universally
measurable subset of a Polish space (in particular, a Borel subset)
and $\mathcal{B}$ the Borel $\sigma$-field on $\Omega$.

Thus, Theorem \ref{mwuv} applies whenever $P$ is perfect,
$\mathcal{A}_n$ is a decreasing sequence of countably generated
$\sigma$-fields, and $\mathcal{A}=\bigcap_n\mathcal{A}_n$. This
particular case, where $\mathcal{A}$ is a tail $\sigma$-field, has
been already proved in Theorem 15 of \cite{BR07}. But Theorem
\ref{mwuv} covers various other real situations.

\begin{ex} {\bf (Increasing unions of tail $\sigma$-fields)} Let $(X_{i,j}:i,j=1,2,\ldots)$ be
an array of real random variables on $(\Omega,\mathcal{B},P)$ and
\begin{equation*}
Z_m^{(n)}=(X_{m,1},\ldots,X_{m,n}),\quad\mathcal{A}_n=\bigcap_m\sigma(Z_m^{(n)},Z_{m+1}^{(n)},\ldots).
\end{equation*}
Then, $\mathcal{A}_1\subset\mathcal{A}_2\subset\ldots$ and, for
each $n$, $\mathcal{A}_n$ is the tail $\sigma$-field of the
sequence $(Z_m^{(n)}:m\geq 1)$. Thus, $\mathcal{A}_n$ meets
condition \eqref{albergo} as far as $\mathcal{B}$ is countably
generated and $P$ perfect. In that case, by Theorem \ref{mwuv},
condition \eqref{2law} holds for
$\mathcal{A}=\sigma(\cup_n\mathcal{A}_n)$.
\end{ex}

\begin{ex}\label{rtrtr} {\bf (Increasing unions of finite groups)} Let
\linebreak $F_1\subset F_2\subset\ldots$ be an increasing sequence
of finite groups of measurable functions of $\Omega$ into itself,
and
\begin{equation*}
\mathcal{A}=\mathcal{A}_{\bigcup_nF_n}=\bigcap_n\mathcal{A}_{F_n}.
\end{equation*}
Suppose $\mathcal{B}$ is countably generated and $P$ perfect. Then,
Theorems \ref{group} and \ref{mwuv} imply that $\mu(\omega)$ is 0-1
on $\mathcal{A}$ for each $\omega\in A_0$, where $A_0\in\mathcal{A}$
and $P(A_0)=1$.
\end{ex}

Finally, as a last and most important example, we mention the {\it
symmetric} $\sigma$-field.

\begin{ex}\label{gxmmj} {\bf (Symmetric $\sigma$-field)} As in Example
\ref{xsetrr}, let
$(\Omega,\mathcal{B})=(\mathcal{X}^\infty,\mathcal{U}^\infty)$
where $(\mathcal{X},\mathcal{U})$ is a measurable space and
$\mathcal{U}$ is countably generated. Denoting $F_n$ the group of
permutations of order $n$, the symmetric $\sigma$-field is
\begin{equation*}
\mathcal{A}=\{B\in\mathcal{B}:f^{-1}B=B\text{ for all }f\in\bigcup_n
F_n\}=\bigcap_n\mathcal{A}_{F_n}.
\end{equation*}
This is just a case of Example \ref{rtrtr}. If $P$ is perfect,
thus, there is $A_0\in\mathcal{A}$ such that $P(A_0)=1$ and
$\mu(\omega)$ is 0-1 on $\mathcal{A}$ for each $\omega\in A_0$.
\end{ex}

\end{document}